\documentclass[12pt]{amsart}
\usepackage{multirow} 
\usepackage{amssymb,amsmath, cite, graphicx, url}
\usepackage[curve]{xypic}
\usepackage{caption}
\makeatletter
\def\@cite#1#2{{\m@th\upshape\bfseries%
[{#1\if@tempswa{\m@th\upshape\mdseries, #2}\fi}]}}
\makeatother
\theoremstyle{plain}
\newtheorem{thm}[subsection]{Theorem}
\newtheorem{cor}[subsection]{Corollary}
\newtheorem{prop}[subsection]{Proposition}
\newtheorem{lem}[subsection]{Lemma}

\theoremstyle{definition}
\newtheorem{rem}[subsection]{Remark}

\newtheorem{war}[subsection]{Warning}

\numberwithin{equation}{subsection}
\newcommand{\bC}{{\mathbb{C}}}
\newcommand{\bQ}{{\mathbb{Q}}}
\newcommand{\bR}{{\mathbb{R}}}

\newcommand{\bH}{{\mathbb{H}}}
\newcommand{\bZ}{{\mathbb{Z}}}

\newcommand{\bN}{{\mathbb{N}}}

\newcommand{\M}{{\mathcal{M}}}

\newcommand{\bk}{{\mathbf{k}}}
\newcommand{\br}{{\mathbf{r}}}

\newcommand{\CP}{{\mathbb C}\!\operatorname{P}^1}
\newcommand{\CPm}{{\mathbb C}\!\operatorname{P}^1\setminus\{z_1, z_2, z_3, z_4\}}

\newcommand{\F}{F} 
\renewcommand{\Im}{\operatorname{Im}}

\begin{document}

\title[Abelian square-tiled surfaces]{Schwarz triangle mappings and Teichm\"uller curves: \\abelian square-tiled surfaces}
%
\author[A.Wright]{Alex~Wright}
\address{Math.\ Dept.\\U. Chicago\\
5734 S. University Avenue\\
Chicago, Illinois 60637}
\email{alexmwright@gmail.com}
%
\date{}

\begin{abstract}

We consider normal covers of $\CP$ with abelian deck group, branched over at most four points. Families of such covers yield arithmetic Teichm\"uller curves, whose period mapping may be described geometrically in terms of Schwarz triangle mappings. These Teichm\"uller curves are generated by abelian square-tiled surfaces.

We compute all individual Lyapunov exponents for abelian square-tiled surfaces, and demonstrate a direct and transparent dependence on the geometry of the period mapping. For this we develop a result of independent interest, which, for certain rank two bundles, expresses Lyapunov exponents in terms of the period mapping. In the case of abelian square-tiled surfaces, the Lyapunov exponents are ratios of areas of hyperbolic triangles. 
\end{abstract}


\maketitle

\setcounter{tocdepth}{1}
\tableofcontents
\newpage

%
%
%
%

\section{Introduction}\label{S:intro}

The general context for this work is the study of rational billiards and translation surfaces. The most well understood examples in this field arise from Teichm\"uller curves, and the most tractable Teichm\"uller curves are the arithmetic ones. The goal of this paper is to study arithmetic Teichm\"uller curves arising from families of abelian covers of $\CP$ branched over four points. These are distinguished by a good geometric description of the period map using Schwarz triangle mappings, and, as we show in \cite{W2}, deep connections to the non-arithmetic Veech-Ward-Bouw-M\"oller Teichm\"uller curves. 

Our study builds upon work by Eskin, Forni, Kontsevich, Matheus, Yoccoz, and Zorich in the cyclic case (details below), and our use of Schwarz triangle mappings to study problems in flat geometry is inspired by work of Bouw and M\"oller \cite{BM}. Our main result is

\begin{thm}
Let $\M$ be any Teichm\"uller curve generated by an abelian square-tiled surface. The Hodge bundle $H^1$ admits a splitting
$$H^1=\bigoplus_{\br} H^1(\br)$$
into flat subbundles $H^1(\br)$ of rank at most 2 with the following properties. 

If $H^{1,0}(\br)$ and $H^{0,1}(\br)$ both have rank one, then the period mapping of $H^1(\br)$ is described by a Schwarz triangle mapping onto a hyperbolic triangle $T$, and the Lyapunov spectrum of $H^1(\br)$ is given by $$\pm\frac{\operatorname{area}(T)}{\pi}.$$
Otherwise $H^1(\br)$ is equal to either its $(1,0)$ or $(0,1)$ part, and its Lyapunov spectrum is all zeros. 
\end{thm}

As far as we know, these Teichm\"uller curves are the only ones in high genus for which the entire Lyapunov spectrum has been computed. The Lyapunov exponents of cyclic square-tiled surfaces were computed using a different method by Eskin-Kontsevich-Zorich \cite{EKZsmall}. Definitions and more precise statements may be found below. 

\textbf{Definitions.} A \emph{Teichm\"uller curve} is an isometrically immersed curve in the moduli space $\M_{g,n}$ of genus $g$ Riemann surfaces with $n$ marked points, with respect to the Teichm\"uller metric. The most basic example is $\M_{0,4}$; the entire moduli space is itself a curve. Any quadratic differential cotangent to a Teichm\"uller curve determines a flat surface which is said to \emph{generate} the Teichm\"uller curve. The cotangent vectors to the point $$(\CP; \{z_1, z_2, z_3, z_4\}) \in \M_{0,4}$$ are multiples of 
\begin{equation}\label{E:q0}
q_0= \frac{(dz)^2}{(z-z_1)(z-z_2)(z-z_3)(z-z_4)}.
\end{equation}
With this flat metric, $\CP$ looks like two isometric parallelograms glued together. For some choice of $z_j$, these parallelograms are squares.

\begin{figure}[h]
\includegraphics[scale=0.2]{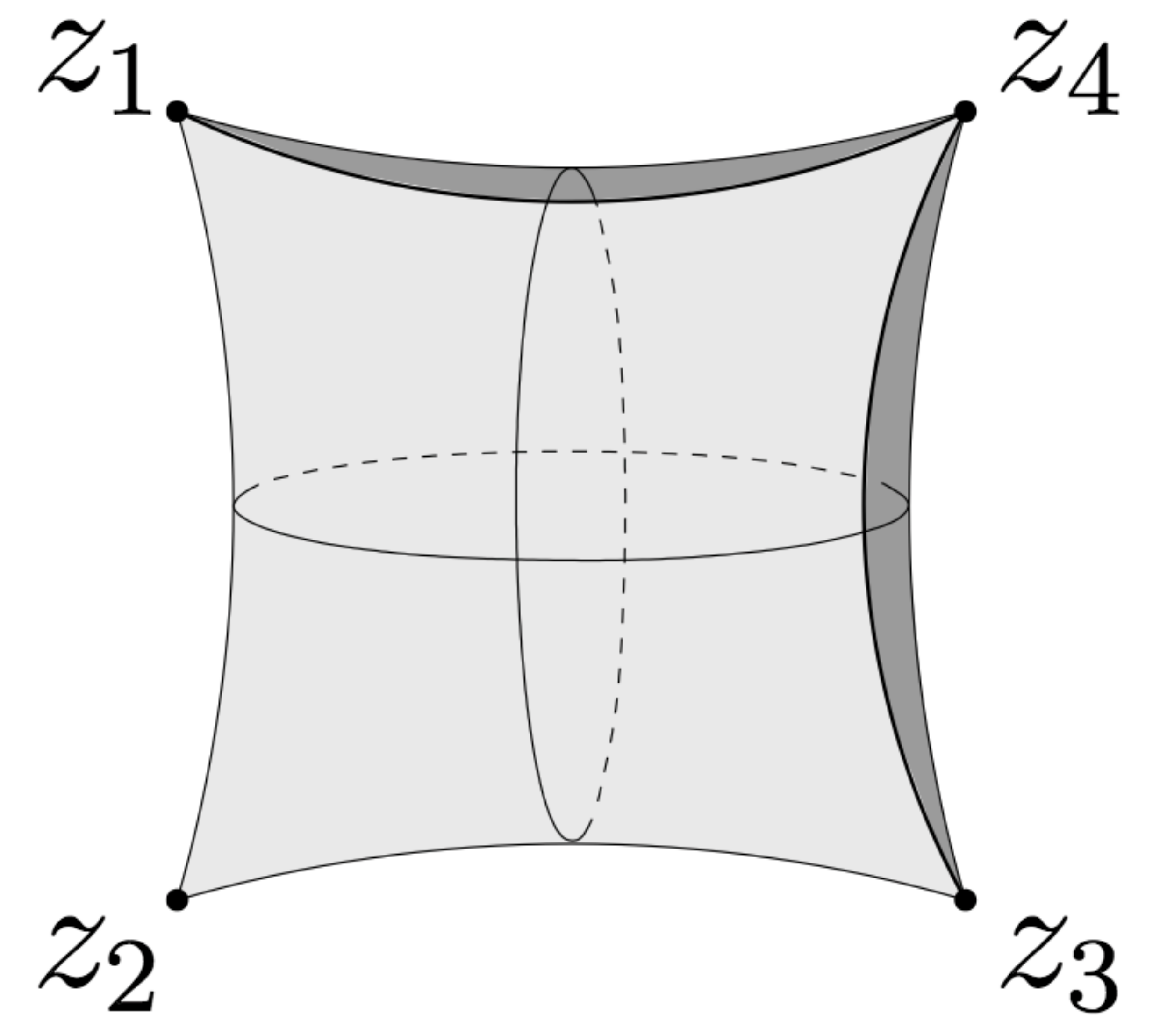}
\caption{The flat pillowcase.} \label{F:FlatPillow}
\end{figure}

An \emph{abelian parallelogram-tiled surface} is a normal cover of $\CP$ which is branched over at most $\{z_1, z_2, z_3, z_4\}$ and has abelian deck group. This cover is equipped with the lift of the flat metric $q_0$ on $\CP$. For some choice of $z_j$, the surface is square-tiled instead of merely parallelogram-tiled. 

Combinatorially, there is no difference between squares and parallelograms, so we will sometimes use ``square-tiled" to mean ``parallelogram-tiled."

Although it is always possible to express an abelian cover as a fibered product of cyclic covers, for us it is crucial to choose more flexible notation which nonetheless remains concrete. Using the Galois correspondence between covers of $\CP$ and field extensions of $\bC(z)$, in Section \ref{S:abcov} we give such notation: given $N>1$, for certain 4 by $m$ matrices $A$ with entries in $\bZ_N=\bZ/(N\bZ)$ we define an abelian square-tiled surface $M_N(A)$. Here $N$, $m$, and $A$ are parameters with no direct geometric interpertation. 

It is possible that $M_N(A)$ is isomorphic to $M_{N'}(A')$ even if $(N,A)\neq (N',A')$. The cover $M_N(A)$ does depend on the branch points $z_j$ (or the $\lambda$ below), but this dependence is supressed in the notation. We have chosen to suppress the dependence to emphasize that the combinatorial structure of the parallelogram tiled-surface remains constant as the branch points vary. 

We may normalize using a M\"obius transformation so that $z_1=0, z_2=1, z_3=\lambda, z_4=\infty$. Varying the parameter $\lambda$ gives a family $\pi:\M_N(A)\to B_0$ of abelian square-tiled surfaces over the base $$B_0=\CP\setminus\{0,1,\infty\}\ni\lambda.$$ The fibers $\pi^{-1}(\lambda)$ of the family $\M_N(A)$ are the abelian covers $M_N(A)$. 

In Section \ref{S:fams} we see that the image of $\M_N(A)$ in moduli space $\M_{g,n}$ is an \emph{arithmetic} Teichm\"uller curve, which we define as a Teichm\"uller curve generated by a square-tiled surface (possibly having non-trivial linear holonomy). Except in unusual cases, we may take $n=0$; typically no marked points are required. 

\textbf{Splitting of the Hodge bundle.} The first cohomology of the fibers of $\pi$ comprise a bundle $H^1$ over $B_0$, commonly known as $R^1\pi_*\bC$. The fiber of $H^1$ over $\lambda\in B_0$ is $H^1(\pi^{-1}(\lambda), \bC)$, where $\pi^{-1}(\lambda)$ is the abelian cover $M_N(A)$ sitting over $\lambda\in B_0$. A similar construction gives a bundle $H^{1,0}$. 

The deck group of $M_N(A)$ is naturally generated by four transformations $T_j$ associated to counterclockwise loops about the $z_j$. We will see that each $T_j$ has order dividing $N$. Let $\xi_N$ be a primitive $N$--th root of unity. For each tuple 

$$\br=(r_1,r_2,r_3,r_4)$$

of integers modulo $N$, we will define the simultaneous eigenspace

$$H^{1,0}(\br)= \bigcap_{j=1}^4  \ker( T_j^*-\xi_{N}^{r_i})\subset H^{1,0}.$$

There is a splitting $$H^{1,0}=\bigoplus_{\br}H^{1,0}(\br).$$ Setting $$H^1(\br)=H^{1,0}(\br)\oplus \overline{H^{1,0}(-\br)},$$ we also have $H^1=\oplus H^1(\br)$. The dimension of $H^1(\br)$ is always at most two. 

\textbf{Period mapping.} The most important case is when both the $(1,0)$ and $(0,1)$ parts of $H^1(\br)$ have dimension 1. Since $H^1(\br)$ is defined topologically in terms of deck transformations, it is a flat subbundle of $H^1$. However, as $\lambda$ varies, the position of $H^{1,0}(\br)$ inside $H^1(\br)$ varies. To record this \emph{variation of Hodge structure} we fix two homology classes $\alpha, \beta\in H_1(M_N(A), \bC)$ and $\omega_\lambda \in H^{1,0}$ (in our case a global section) and compute the period mapping

$$\tau(\lambda)=\frac{ \int_{\alpha} \omega_\lambda}{ \int_{\beta}\omega_\lambda}.$$

The period mapping is a holomorphic map 

$$\tau: \CP\setminus\{0,1,\infty\}\to \CP.$$

 Changing $\alpha$ and $\beta$ changes $\tau$ by postcomposition with a M\"obius transformation. Note that parallel translating $\alpha$ and $\beta$ around a loop in $\CP\setminus\{0,1,\infty\}$ results in monodromy, so $\tau$ is not globally defined on $\CP\setminus\{0,1,\infty\}$. But if we restrict ourselves to the simplify connected set $\bH \subset\CP\setminus\{0,1,\infty\},$ there are no issues with monodromy, so $\tau$ is well defined (up to M\"obius transformations) on $\bH$.  
 
In Sections \ref{S:HGDE} and \ref{S:periodmap} we recall that $\tau$ is a \emph{Schwarz triangle mapping}, in particular a biholomorphism from $\bH \subset \CP\setminus\{0,1,\infty\}$ to a hyperbolic triangle in $\bH$. 

All the information about a family of Riemann surfaces is contained in its period mapping, and this is one of very few situations where such an explicit and elegant description is possible. In particular, this description will allow us to compute Lyapunov exponents. 

\textbf{Lyapunov exponents.} Given a generic cohomology class in $H^1(\br)$, we may parallel transport along longer and longer geodesics on the Teichm\"{u}ller curve corresponding to $\M_N(A)$. The exponential rate of growth of the cohomology class (with respect to the Hodge norm) is the top \emph{Lyapunov exponent} $\ell$ of $H^1(\br)$. 

General results applied to this situation relate $\ell$ to the degree of the line bundle $H^{1,0}(\br)$. These results were developed by Forni and Bouw-M\"oller following numerical experiments of Zorich and a breakthrough observation of Kontsevich \cite{K, F, BM}. 

In all of dynamics this is one of only very few examples where Lyapunov exponents may be explicitly computed. In flat geometry the significance of Lyapunov exponents is twofold: they describe both the dynamics of the Teichm\"uller geodesic flow on moduli space, and the deviation of ergodic averages for straight line flow on the translation surface \cite{F}.

In Section \ref{S:lyaps}, we prove the following variant of the above results. In the theorem below, $W=H^1(\br)$ is the example we have in mind. 

\begin{thm}\label{T:periodlyap}
Let $\M$ be a Teichm\"uller curve, and let $W\subset H^1$ be a flat subbundle with $\dim_\bC W^{1,0}=1=\dim_{\bC} W^{0,1}$ and $\dim_\bC W=2$. The nonnegative Lyapunov exponent $\ell$ of $W$ is 
$$\ell=\frac1{Area(\M)} \int_\M \| \tau'(\lambda)\|^2_{hyp} dg_{hyp}(\lambda),$$
where $\tau$ is the period mapping of $W$. 
\end{thm}

Because the period mapping is well defined up to post-composition with M\"obius transformations, the hyperbolic norm $\| \tau'(\lambda)\|_{hyp}$ of the derivative is well defined. 

For abelian square-tiled surfaces, $\M=\CP\setminus\{0,1,\infty\}$ is equipped with a hyperbolic metric via uniformization, and $dg_{hyp}$ is the corresponding hyperbolic volume form.

\begin{thm}\label{T:arealyap}
When the $(1,0)$ and $(0,1)$ parts of $H^1(\br)$ each have dimension one, then the nonnegative Lyapunov exponent is equal to $$\ell=\frac{A}{\pi},$$ where $A$ is the hyperbolic area of the Schwarz triangle $\tau(\bH)$. Otherwise, the Lyapunov exponents of $H^1(\br)$ are zero. 
\end{thm}

\begin{proof}[\textbf{\emph{Proof sketch.}}]
First consider the case where the $(1,0)$ and $(0,1)$ parts of $H^1(\br)$ each have dimension one. The period map $\tau$ of $H^1(\br)$ maps $\bH\subset \CP\setminus\{0,1,\infty\}$ to a hyperbolic triangle of area $A$. The hyperbolic manifold $\CP\setminus\{0,1,\infty\}$ is composed of two ideal triangles of area $\pi$, which are the upper and lower half-planes in $\CP\setminus\{0,1,\infty\}$. The change of variables formula gives that 
$$\int_{\bH\subset \CP\setminus\{0,1,\infty\}} \|\tau'(\lambda)\|^2_{hyp} dg_{hyp}(\lambda)$$ is the hyperbolic area $A$ of the triangle $\tau(\bH)$, from which Theorem \ref{T:periodlyap} gives the result.

In the second case, $H^1(\br)$ is equal to either its $(1,0)$ or $(0,1)$ part. It follows that the Hodge bilinear form on $H^1(\br)$ is positive definite, so the monodromy is contained in a compact group. 
\end{proof}

In Section \ref{S:lyaps} we give an algorithmic restatement of Theorem \ref{T:arealyap} which computes all the Lyapunov exponents for any abelian square-tiled surface. 

\textbf{Holonomy double cover.} When the quadratic differential defining the flat structure is not the square of an abelian differential, the \emph{holonomy double cover} is the smallest cover where the lift of the quadratic differential is a square. In Section \ref{S:flat}, we show that the holonomy double cover of an abelian square-tiled surface is again an abelian square-tiled surface. 

\textbf{Combinatorics.} In Section \ref{S:flat} we give a combinatorial model of $M_N(A)$. The main point is that the square tiles may be labeled by elements of the deck group, which is isomorphic to the column span of $A$. 

\textbf{Hypergeometric differential equations.} To show that the period mapping of a $H^1(\br)$ is a Schwarz triangle mapping, we consider a global section $\omega_\lambda$ of $H^{1,0}(\br)$. When $H^1(\br)$ is two dimensional, $\omega_\lambda$ and its first two derivatives $\omega_\lambda', \omega_\lambda''$ must be linearly dependent, and hence satisfy some linear relation which depends on $\lambda$. That is, $\omega_\lambda$ satisfies some second order differential equation, which turns out to be the HGDE. The periods 

$$ \int_{\alpha} \omega_\lambda, \quad \int_{\beta}\omega_\lambda$$

of $\omega_\lambda$ also satisfy this HGDE. 

The period mapping is thus the ratio of two solutions to a HGDE, which is by definition a Schwarz triangle mapping. In Section \ref{S:HGDE} we recall the proof that a Schwarz triangle mapping maps the upper half-plane biholomorphically to a triangle.  

\textbf{Context.} In \cite{F}, Forni showed that generically, the Lyapunov exponents of a translation surface are all nonzero. Interest in cyclic square-tiled surfaces grew out of a small number of examples due to Forni, Matheus and Zorich where the Lyapunov spectrum is \emph{completely degenerate}, that is, has as many zero exponents as possible \cite{Fex1, FM, FMZ}. The monodromy of some of these examples has been studied by Matheus and Yoccoz \cite{MY}. 

It is expected that the only Teichm\"uller curves (generated by abelian differentials) with completely degenerate Lyapunov spectrum are the known examples arising from cyclic square-tiled surfaces. M\"oller has proved this in every genus except $g=5$ \cite{M5}.

As an application of their longstanding project on Lyapunov exponents, Eskin-Kontsevich-Zorich computed the Lyapunov exponents for cyclic square-tiled surfaces \cite{EKZbig, EKZsmall}. Using the theory of abelian square-tiled surfaces which we develop, their methods could be used instead of Theorem \ref{T:periodlyap} to compute the Lyapunov spectrum. Their methods are more algebro-geometric and
clarify the role of the Deligne extension, whereas our methods are more visually geometric and clarify the role of the period mapping. 

Arithmetic Teichm\"uller curves have been extensively studied, but usually under the assumption that $n=0$ (no marked points) and trivial linear holonomy. The terminology \emph{arithmetic} is apt because such curves are uniformized by subgroups of $PSL(2,\bZ)$ \cite{GJ}; at least under the usual assumptions they admit an action of the absolute Galois group of $\bQ$ \cite{M4}; and, in analogy with a result of Margulis for lattices, they are the only Teichm\"uller curves whose commensurability classes do not admit minimal representatives \cite{M4}.

\textbf{References} For an introduction to translation surfaces and rational billiards see, for example, the surveys \cite{MT,Z}.

The theory of cyclic and abelian covers of $\CP$ and the connection to HGDE and Schwarz triangle mappings is classical. We include foundational results on these topics because our notation is new, because we have not found suitable references, and because we do not expect our audience to be familiar with these results. 

\textbf{Terminology.} Frequently, definitions of square-tiled surface require trivial linear holonomy. In contrast, following \cite{FMZ, EKZsmall}, our square-tiled surfaces may be generated by quadratic differentials which are not squares of abelian differentials.

\textbf{Acknowledgements.} This research was supported in part by the National Science and Engineering Research Council of Canada, and was partially conducted during the Hausdorff Institute's trimester program ``Geometry and dynamics of Teichm\"uller space." The author thanks the Hausdorff Institute for its hospitality. 

This paper grew out of the problem of computing the Lyapunov exponents of the holonomy double cover of a cyclic square-tiled surface, suggested to the author by his thesis advisor Alex Eskin. The author thanks Alex Eskin for this suggestion, and for his guidance and support. Special thanks are also due to Howard Masur, Martin M\"oller, and Anton Zorich for their instruction and encouragement. The author is grateful to Matt Bainbridge, Irene Bouw, Madhav Nori and Daniel Studenmund for helpful and interesting discussions.

The author thanks Anton Zorich for Figure \ref{F:FlatPillow} and Jennifer Wilson for Figure \ref{F:triangleslyap}.


\section{Abelian covers}\label{S:abcov}

In this section we give basic results on abelian covers of $\CP$ branched over at most four points.


\subsection{Parameterization.}\label{SS:param} Let $N>0$ and let $A$ be an $m$ by $4$ matrix with entries $a_{i,j}$ in $\bZ_N=\bZ/(N\bZ)$. We will assume that the sum of the four columns of $A$ is zero. That is, we require that the sum of the four entries in any row is zero. 

Let $\overline{\bC(z)}$ denote the algebraic closure of $\bC(z)$. For each $i=1, \ldots, m$, pick a function $w_i\in \overline{\bC(z)}$ with 
\begin{equation}\label{E:wi}
w_i^N=\prod_{j=1}^4 (z-z_j)^{\tilde{a}_{i,j}},
\end{equation}
 where $\tilde{a}_{i,j}$ is the lift of $a_{i,j}$ to $\bZ\cap [0,N)$. 

Let $M_N(A)$ be the smooth nonsingular algebraic curve whose function field is $\bC(z)[w_1,\ldots, w_m]$. This Riemann surface comes equipped with a branched covering map to the $\CP$ with function field $\bC(z)$, and this cover is normal with abelian deck group. The dependance on the set of branch points $\{z_1,\ldots, z_4\}$ is suppressed in the notation $M_N(A)$.

Conversely, consider a Riemann surface $X$ which is a normal cover of $\CP$, branched over $\{z_1,z_2, z_3, z_4\}$, with abelian deck group. Its function field is an abelian extension of $\bC(z)$, the function field of $\CP$.

A basic result in Galois theory gives that a degree $D$ cyclic extension of a field which contains all $D$-th roots of unity is obtained by taking a $D$-th root. Taking into account the branch points, this result gives that $X=M_N(A)$ for some choice of $N$ and $A$ as above. 

We emphasize that when using the notation ``$M_N(A)$" we will always require that the sum of the columns rows of $A$ is zero. This is to ensure that $M_N(A)$ is not branched over $z=\infty$ in addition to the four $z_j$. Furthermore, $z$ will always denote a parameter of the $\CP$ of which $M_N(A)$ is an abelian cover, and $w_i$, for $i=1,\ldots, m$, will always denote a fixed element of the function field as in equation \ref{E:wi}. 

Later we will endow $M_N(A)$ with the structure of a flat surface, sometimes with finitely many unlabeled marked points. (In Section \ref{SS:stratum} will see that if $A$ has no zero columns, then no marked points are required.)


\subsection{Galois group.}\label{SS:galgroup} Let $F$ be the function field of $M_N(A)$, considered as a Galois extension of the field $\bC(z)$. To study this extension, it is helpful to consider a larger field $E\supset F$, whose Galois theory will be particularly clear. We will find a basis for $F$ over $\bC(z)$ as a subset of a basis for $E$ over $\bC(z)$, and we will compute the Galois group $Gal_{\bC(z)}(F)$ as the image of $Gal_{\bC(z)}(E)$ under the restriction map. 

Consider then the field $E=\bC(z)[u_1, u_2, u_3, u_4]$ where $u_j$ is an $N$-th root of $z-z_j$ for each $j=1,2,3,4$. Note that $F$ is contained in $E$. 

The set of products $\prod_{j=1}^4 u_j^{c_j}$ with $0\leq c_j<N$ form a basis for $E$ over $\bC(z)$. Let $\xi_N$ be a primitive $N$-th root of unity. The Galois group of $E$ over $\bC(z)$ is $\bZ_N^4$, with generators $\psi_j$, $j=1,2,3,4$, where $\psi_j(u_k)=u_k$ if $j\neq k$ and $\psi_j(u_j)=\xi_N u_j$.

Consider now the collection $B\subset F$ consisting of $\prod_{j=1}^4 u_j^{\tilde{r_j}}$ for each $(r_1, r_2,r_3, r_4)$ in the row span of $A$. This collection $B$ spans $F$, and is automatically linearly independent since it is a subset of a basis for $E$ over $\bC(z)$. Thus, $B$ is a basis for $F$ over $\bC(z)$. It follows that the degree of the cover $M_N(A)$ is the size of the row span of $A$. 

There is a surjective group homomorphism $Gal_{\bC(z)}(E)\to Gal_{\bC(z)}(F)$ given by restriction of field automorphisms from $E$ to $F$.

To compute the restriction of $\psi=\prod_{j=1}^4 \psi_j^{p_j}$ to $F$, it suffices to compute its action on each $w_i$. Note that $\psi$ acts on $w_i$ by multiplication with $\xi_N^{\sum_{j=1}^4 a_{i,j}p_j}$. Consider the column vector $(c_1, \ldots, c_m)^t\subset \bZ_N^m$ given as the sum of $p_j$ times the $j$-th column of $A$, $j=1,2,3,4$. This column vector records the action of $\psi$ on each $w_i$; we see that $\psi$ acts on $w_i$ as multiplication by $\xi_N^{c_i}$. 

It follows that $Gal_{\bC(z)}(F)$ is isomorphic to the column span of $A$.

Since $F$ is a Galois extension of $\bC(z)$, the size of the Galois group is equal to the degree of the extension.  That is,  the size of the row span of $A$ is equal to the size of the column span of $A$. This fact is outside the realm of linear algebra, since we are working modulo $N$, but it can also seen from an elementary argument using Smith normal form.


\subsection{Partial order on the set of abelian covers.}\label{SS:order} The cover $M_N(A)$ is covered by $M_{N'}(A')$, in a way compatible with the covering maps to $\CP$, exactly when the function field of $M_N(A)$ is contained in that of $M_{N'}(A')$.  This occurs if and only if there exist $k', k\in \bN$  such that $k'N'=kN$ and the row span of $kA$ is contained in the row span of $k'A'$. (If $A$ is a matrix with entries in $\bZ_N$, then $kA$ may naturally be considered as a matrix with entries in $\bZ_{kN}$.) Two abelian covers $M_N(A)$ and $M_{N'}(A')$ are isomorphic if and only if  there exist $k', k\in \bN$  such that $kN=k'N'$ and $kA$ and $k'A'$ have the same row span. 

An expression ``$M_N(A)$" is called a \emph{presentation} of the abelian cover it determines. Every abelian cover admits many presentations. It is important for applications in \cite{W2} that we not put any restrictions on presentations, even though every abelian cover admits a presentation of a quite restricted form.


\subsection{Genus.}\label{SS:genus} The deck group of $M_N(A)$ is isomorphic to the Galois group of the field extension. The Galois group is naturally isomorphic to the column span of $A$, which has the same size as the row span of $A$ (Section \ref{SS:galgroup}). The degree $d$ of the cover $M_N(A)\to\CP$ is thus size of the column span of $A$. The order of ramification over $z_j$ is the order of the loop around $z_j$ in the deck group of $M_N(A)$. The image of the loop around $z_j$ in the deck group is the transpose of $(a_{1,j}, \ldots, a_{m,j})$, and  hence has order $N/\gcd(N,a_{1,j}, \ldots, a_{m,j})$. (Technically, lifts of the $a_{i,j}$ to $\bZ$ should be used in the $\gcd$, but the result does not depend on the choice of lifts.) Now, from the Riemann-Hurwitz formula
we compute
\begin{eqnarray}
g
&=&
1+ d \left(1 - \frac1{2N} \sum_{j=1}^4 \gcd(N, a_{1,j}, \ldots, a_{m,j})\right). \label{E:genus}
\end{eqnarray}


\subsection{A basis of holomorphic one-forms.}\label{SS:basis} 
Define, for $\br=(r_1, r_2, r_3, r_4)$ in the row span of $A$, $$t_j(\br)=\left\{ \frac{\tilde{r}_j}{N} \right\},$$ where $\tilde{r}_j$ is a lift of $r_j$ to $\bZ$, and $\{\cdot\}$ denotes fractional part. The fractional part $\{ x\}$ is equal to $x$ minus the largest integer less than or equal to $x$. Define
$$t(\br) = t_1(\br)+t_2(\br)+t_3(\br)+t_4(\br). $$

Since the sum of the entries in any row of $A$ is zero, we have that $t(\br)\in \{1,2,3\}$ if $\br\neq 0$. 

\begin{lem}\label{L:omegahol} 
Let $\br$ in the row span of $A$ be nonzero. Set 
$$\omega=(z-z_1)^{-t_1(\br)}(z-z_2)^{-t_2(\br)}(z-z_3)^{-t_3(\br)}(z-z_4)^{-t_4(\br)}dz.$$
The meromorphic-one form  $p(z)\omega$ with  $p(z)\in \bC(z)$ is a holomorphic one-form on $M_N(A)$ if and only if $p$ is a polynomial of degree at most $t(\br)-2$. 
\end{lem}

\begin{proof}
Write $p(z)=p_0(z) \prod_{j=1}^4 (z-z_j)^{s_j}$, where $p_0$ does not have any roots or poles at any of the $z_j$. Assume $p(z) \omega$ is holomorphic. So $p_0$ is a polynomial. 

Let $m_j$ be the order of ramification of $M_N(A)$ over $z_j$. We may choose a local parameter $u$ near a lift of $z_j$ so that $u^{m_j}=(z-z_j)$. Near this lift of $z_j$ we compute that $p(z)\omega$ is proportional to $$u^{(s_j-t_j(\br)+1)m_j -1} du.$$ Hence $p(z)\omega$ is holomorphic at lifts of $z_j$ if 
\begin{equation}\label{E:holz}
s_j\geq \frac1{m_j} +t_j(\br)-1.
\end{equation}
Recall $m_j=N/\gcd(N,a_{1,j}, \ldots, a_{m,j})$, and note that $t_j(\br)<1$ is a multiple of of $1/m_j$, so 
$$0<\frac{1}{m_j}+ t_j(\br)\leq 1.$$
Since $s_j$ is an integer, \ref{E:holz} is equivalent to $s_j\geq 0.$

Hence $p$ is a polynomial, say of degree $d$. The condition that $p(z)\omega$ is holomorphic over lifts of $\infty$ is $d\leq t(\br)-2$, which gives the result.
\end{proof}

\textbf{Splitting of $H^{1,0}$.} Set 
$$H^1=H^1(M_N(A), \bC),\quad H^{1,0}=H^{1,0}(M_N(A)),\quad H^{0,1}=H^{0,1}(M_N(A)).$$ 
Let $T_j, j=1,2,3,4$ be the deck transformation of $M_N(A)$ corresponding to the counterclockwise loop about $z_j$. So $\prod_{j=1}^4 T_j=1$, and $T_j(w_i)=\xi_N^{a_{i,j}} w_i$. For $\br=(r_1,r_2,r_3, r_4)\in \bZ_N^4$, we may define 
$$H^1(\br)= \bigcap_{j=1}^4  \ker( T_j^*-\xi_{N}^{r_i})\subset H^1,$$
where $T_j^*: H^1\to H^1$ is the induced action on cohomology. The maps $T_j$ preserve the complex structure, and hence induce maps $T_j^*: H^{1,0}\to H^{1,0}$ on the space of holomorphic one-forms. This action preserves the Hodge bilinear form, 
\begin{equation}\label{E:hodge}
\langle \omega_1, \omega_2\rangle = \frac{i}2 \int \omega_1\wedge \overline{\omega_2},
\end{equation}
which is positive definite on $H^{1,0}$ and negative definite on $H^{0,1}$. We may conclude that the commuting linear maps $T_j^*$ on $H^{1,0}$ are diagonalizable and hence simultaneously diagonalizable. If we set 
$$H^{1,0}(\br)= \bigcap_{j=1}^4  \ker( T_j^*-\xi_{N}^{r_i})\subset H^{1,0},$$
we find that $$H^{1,0}=\bigoplus_{\br\in \bZ_N^4} H^{1,0}(\br).$$ In fact, it is not hard to see that the summation may be restricted to $\br$ in the row span of $A$. Furthermore, $H^{1,0}(0)$ is trivial, since any homomorphic one-form invariant under all $T_j^*$ descends to a holomorphic one form on $\CP$, and there are no holomorphic one forms on $\CP$.

Since $H^{1,0}=\bigoplus_{\br} H^{1,0}(\br)$ it follows that $H^1=\bigoplus_{\br} H^1(\br)$. The summations are over nonzero $\br$ in the row span of $A$.  We may also define $H^{0,1}(\br)= \bigcap_{j=1}^4\ker (T_j^*-\xi_{N}^{r_i})\subset H^{0,1}$ using the induced action of the $T_j$ on $H^{0,1}$, and we similarly have $H^{0,1}=\bigoplus_{\br} H^{0,1}(\br)$.

Lemma \ref{L:omegahol} gives: 

\begin{prop}\label{P:dim}
For $\br$ in the row span of $A$ and not zero, we have 
$$\dim_\bC H^{1,0}(\br)=t(-\br)-1 \in \{0,1,2\},$$
and
$$\dim_\bC H(\br)=t(\br)+t(-\br)-2\in \{0,1,2\}.$$
If $t(-\br)\geq 2$, then $$\omega= \prod_{j=1}^4(z-z_j)^{-t_j(-\br)} dz \in H^{1,0}(\br),$$ and if $t(-\br)=3$, then $z \omega \in H^{1,0}(\br)$ as well. 
\end{prop}

\begin{proof}
$T_j$ also induces an action on the function field $F$ of $M_N(A)$. (In the notation of Section \ref{SS:galgroup}, the action of $T_j$ on $F$ is $\psi_j$.) This action is linear over $\bC(z)$, and $\prod_{j=1}^4(z-z_j)^{-t_j(-\br)}$ spans the simultaneous eigenspace where $T_j$ acts by multiplication by $\xi_N^{r_i}$, $j=1,2,3,4$. Hence any $\eta\in H^{1,0}(\br)$ is of the form $$\eta=p(z) \prod_{j=1}^4(z-z_j)^{-t_j(-\br)} dz,$$ with $p(z)\in \bC(z)$. Lemma \ref{L:omegahol} gives that such an $\eta$ is holomorphic if and only if $p(z)$ is a polynomial of degree at most $t(-\br)-2$. 

Since $H^{0,1}(\br)=\overline {H^{1,0}(-\br)}$ and $H^1(\br)= H^{1,0}(\br)\oplus H^{0,1}(\br)$, the formula for $\dim_\bC H(\br)$ follows from the formula  for $\dim_\bC H^{0,1}(\br)$. 
\end{proof}


\section{The flat structure on abelian covers}\label{S:flat} 

In this section we describe the flat square-tiled structure on abelian covers $M_N(A)$.

Recall the flat pillowcase metric on $\CP$ is given by the quadratic differential  
\begin{equation}\label{E:q0}
q_0= \frac{(dz)^2}{(z-z_1)(z-z_2)(z-z_3)(z-z_4)}.
\end{equation}

Equipping the abelian cover $M_N(A)$ with the lift $q$ of the quadratic differential $q_0,$ 
we obtain an abelian square tiled-surface. Later we will consider $M_N(A)$ as a Riemann surface with not only a flat structure but sometimes also with finitely many unlabeled marked points (the poles of $q$).


\subsection{Stratum.}\label{SS:stratum} Let $d$ be the size of the row span of $A$. Recall that the order of ramification above $z_j$ is $N/\gcd(N,a_{1,j}, \ldots, a_{m,j})$. Hence, $M_N(A)$ lies in the stratum of half translation surfaces with $$\frac{\gcd(N,a_{1,j}, \ldots, a_{m,j})}{N} d$$ singularities of  cone angle $$ \frac{N}{\gcd(N,a_{1,j}, \ldots, a_{m,j})}\pi$$ for each $j=1,2,3,4$. 

In particular, notice that $M_N(A)$ has cone angles of $\pi$ if and only if $A$ has a column of zeros.


\subsection{Linear holonomy and the holonomy double cover.}\label{SS:holonomy} The holonomy double cover of a half-translation surface $(X,q)$ is defined to be the smallest cover of $(X,q)$ where $q$ lifts to the square of an abelian differential. This is exactly the smallest cover with trivial linear holonomy, and is at most a two-fold cover of $(X,q)$.  

Consider now the differential $q$ on $M_N(A)$. Since $dz$ is globally defined, $q$ is the square of an abelian differential if and only if the function $(z-z_1)(z-z_2)(z-z_3)(z-z_4)$ has a global square root on $M_N(A)$. This happens if and only if $M_N(A)$ covers $M_2(1,1,1,1)$. 

Define $N^+$ and  $A^+$ as follows. If $N$ is even, then $N^{+}=N$ and $A^+$ is obtained by adding a row of $N/2$'s to $A$. If $N$ is odd, $N^+=2N$, and $A^+$ is obtained by adding a row of $N$'s to $2A$. 

\begin{prop}
The holonomy cover of $M_N(A)$ is $M_{N^+}(A^+)$.
\end{prop}

The proof uses the partial order on the set of abelian covers. 

\begin{proof}
$M_{N^+}(A^+)$ is equal to $M_N(A)$ if and only if $M_N(A)$ covers $M_2(1,1,1,1)$, that is, if and only if $q$ is the square of an abelian differential on $M_N(A)$. 

Otherwise, $M_{N^+}(A^+)$ is a twofold cover of $M_N(A)$ to which $q$ lifts to a square of an abelian differential. 
\end{proof}

\begin{cor}
$M_N(A)$ has non-trivial linear holonomy if and only if $N$ is even and $(N/2, N/2, N/2, N/2)$ is in the row span of $A$. 
\end{cor}


\subsection{A combinatorial model.}\label{SS:comb} This section is inspired by \cite{FMZ}, where a related combinatorial model is given for cyclic square-tiled surfaces.

The flat pillowcase $(\CPm, q_0)$ may be (non-canonically) divided into a front and back, so that the front and back are isometric parallelograms. Fix a choice of front and back. The front of the flat pillowcase $(\CPm, q_0)$ is declared white, and the back black. Fixed lifts, by convention with the black square (parallelogram) bordering the white along a lift of the edge joining $z_3$ to $z_4$, of the two squares (parallelograms) of $(\CPm, q_0)$ will act as a sort of fundamental domain for $M_N(A)$, and will both be labelled by zero. Other squares will be labelled by elements in the column span of $A$ with a subscript of $w$ or $b$ to indicate color. 

Consider a vector $\textbf{q}=\left( q_1 , q_2 , \ldots, q_m\right)^t$ in the column span of $A$. This can be naturally considered as a deck transformation. 

 The square $\textbf{q}_w$ is defined as the image of the fundamental white square (labelled by zero) under this deck transformation. The black square $\textbf{q}_b$ is similarly defined.
 

In this way all the squares of $M_N(A)$ are labeled. Note that the total number each of white and black squares is the size of the column span of $A$. Furthermore the action of the deck group, which is the column span of $A$, is given by simply adding an element of the column span of $A$ to each label while leaving the subscripts unchanged.


\section{Families of abelian covers}\label{S:fams}

A family of Riemann surfaces over a base $B_0$ is a holomorphic submersion $\pi:\M\to B_0$ whose fibers are Riemann surfaces. 

In this section we consider a family of curves $\M_N(A)\to \CP\setminus\{0,1,\infty\}$ whose fibers are isomorphic to $M_N(A)$ for some choice of $z_j$. More specifically, consider
\begin{equation}\label{E:fam} 
\begin{split}
\{& (\lambda, z, w_1, \ldots, w_m) \in \CP\setminus\{0,1,\infty\} \times \bC \times \bC^m: \\ 
&w_i^N= z^{\tilde{a}_{1,j}} (z-1)^{\tilde{a}_{2,j}} (z-\lambda)^{\tilde{a}_{3,j}}, \quad j=1\ldots m \quad\}.
\end{split}
\end{equation}
This is a family of possibly singular plane curves, over the base $B_0= \CP\setminus\{0,1,\infty\}$. The map $\pi$ is given by $\pi (\lambda, z, w_1, \ldots, w_m)=\lambda$. We define $\M_N(A)$ to be the family given by taking the normalization of the projectivization of fibers and, on each fiber (Riemann surface), marking the preimages of $z=0,1,\lambda, \infty$ where no branching occurs. 

There is a vector bundle over $\M_N(A)$ with fiber $H^1(X, \bC)$ over $X\in \M_N(A)$.  This bundle is typically called $R^1 \pi_* \bC$, but we will denote it $H^1$. In this section we will also describe a splitting of this bundle. 

\begin{war}
Because the moduli space $\M_{g}$ is not a fine moduli space, neither $\M_{g}$ nor any subset of $\M_{g}$ is a family of curves. For this reason, without using definitions adapted to stacks, the ``bundle over $\M_{g}$ whose fiber over a Riemann surface is its first cohomology" is in fact not a well defined vector bundle. These problems are very mild, and disappear after replacing $\M_{g}$ with an appropriate finite cover. 

Nonetheless, a family of curves should not be confused with its image in moduli space. Given a family, the bundle $R^1 \pi_* \bC$ over $B$, whose fiber over $b\in B$ is $H^1(\pi^{-1}(b), \bC)$, is rigorously defined. Moreover, this bundle endowed with the Gauss-Manin connection is a flat bundle. However, there may be several families over a given base $B$ which give the same map $B\to \M_{g}$ but for which the bundles $R^1 \pi_* \bC$ are not isomorphic as flat bundles, that is, have non-isomorphic monodromy. Such families are fiberwise isomorphic but non-isomorphic. 
\end{war}


\subsection{$\M_N(A)$ as an arithmetic Teichm\"uler curve.}\label{SS:arith} When $M_N(A)$ has cone angles of $\pi$, these occur at simple poles of the quadratic differential $q$ defining the flat metric. We wish for $q$ to be a cotangent vector of the moduli space containing $M_N(A)$. 

Cotangent vectors to $\M_{g,n}$ are holomorphic one forms with at most simple poles at the marked points. This is why, on each Riemann surface in the family $\M_N(A)$, we have marked all poles of $q$, equivalently all points with cone angle $\pi$. Whenever $A$ has no zero columns there are no cone angles of $\pi$, and hence no points need to be marked (Section \ref{SS:stratum}). Nonetheless we emphasize that, when marked points are required, $\M_N(A)$ is considered to be a family in the moduli space $\M_{g,n}$ of genus $g$ Riemann surfaces with $n>0$ unlabeled marked points.  

If $g\in SL_2(\bR)$, and $X$ is a fiber (Riemann surface) of $\M_N(A)$, and $q$ the lift of $q_0$ to $X$, then $g\cdot (X, q)$ is again a cover of $\CP$ with the same monodromy. So $g\cdot (X, q)$ is isomorphic to some fiber of $\M_N(A)$. Hence the image of $\M_N(A)$ in $\M_{g,n}$ is the projection of the closed $SL_2(\bR)$ orbit of a square-tiled surface. In other words, the image of $\M_N(A)$ in $\M_{g,n}$ is an arithmetic Teichm\"uller curve. 

\begin{rem}\label{R:3ps}
The base of the family $\M_N(A)$ is the three times punctured sphere, which has a unique hyperbolic structure. The image of $\M_N(A)$ in moduli space is an isometrically immersed hyperbolic three times punctured sphere. The immersion need not be generically one-to-one, but it is at most six-to-one, because the three times punctured sphere is uniformed by the level two congruence subgroup, which has index six in $PSL_2(\bZ)$.  

In other words, the projective affine group of an abelian square-tiled surface sits in between the level two congruence subgroup and $PSL_2(\bZ)$. 
\end{rem}


\subsection{Splitting of $H^1$.}\label{SS:split2} By abuse of notation, we denote the bundles over the base of $\M_N(A)$ with fibers $H^1(\br)$, $H^{1,0}(\br)$ and $H^{0,1}(\br)$ (defined in Section \ref{SS:basis}) by the same notation $H^1(\br)$, $H^{1,0}(\br)$ and $H^{0,1}(\br)$. Lemma \ref{L:omegahol} gives the rank of these bundles. All of the splittings given in Section \ref{SS:basis}, for example $H^1=\oplus_{\br} H^1(\br)$, are valid for the corresponding bundles.

The bundle $H^1(\br)$ is a \emph{complex VHS}, which for a subbundle of $H^1$ simply means that it decomposes into its $(1,0)$ and $(0,1)$ parts, but that there need not be any underlying real vector bundle. In particular, note that the $(1,0)$ and $(0,1)$ parts of $H^1(\br)$ are usually not conjugate, since $\overline{H^{0,1}(\br)}= H^{1,0}(-\br)$.

The subbundle $H^1(\br)$ is covariantly constant with respect to the Gauss-Manin connection. That is, it is a flat subbundle of $H^1$. 

The goal of the next two sections is to define and describe the period mapping of each $H^1(\br)$.


\section{Schwarz triangle mappings}\label{S:HGDE}

In this section we review the theory of hypergeometric differential equations, focussing on Schwarz triangle mappings. A standard reference on this topic is \cite{Y}; we also recommend the unpublished notes by Frits Beukers \cite{B, B2}, and, for general background on Schwarz-Christoffel mappings, the book \cite{DT}. In the interest of brevity we have not included the general theory of Fuchsian differential equations. Fuchsian differential equations are linear differential equations with only regular singular points, and HGDE are the standard form of second order Fuchsian differential equations with three singular points.

\subsection{Hypergeometric functions.}\label{SS:HGfunctions} For $a,b,c\in \bR$ with $c$ not zero and not a negative integer, the Gauss hypergeometric function $\F(a,b;c; \lambda)$ is defined by 
$$\F(a,b;c; \lambda) = \sum_{n=0}^\infty \frac{(a)_n (b)_n}{(c)_n} \frac{\lambda ^n}{n!}.$$
The Pochhammer symbol is defined by $(x)_n=x(x+1)\cdots (x+n-1)$ with the convention $(x)_0=1$. The radius of convergence of $\F(a,b;c; \lambda)$ is 1 unless $a$ or $b$ is a nonpositive integer, in which case it is $\infty$. Throughout this section we will assume that $c$ is not zero and not a negative integer. 

The function $\F(a,b;c; \lambda)$ satisfies the HGDE
\begin{equation}\label{E:DE}
\lambda(\lambda-1)\frac{d^2}{d\lambda^2}+\left[(a+b+1)\lambda-c\right]\frac{d}{d\lambda} + ab = 0.
\end{equation}
The function $$\lambda^{1-c}\F(a+1-c, b+1-c; 2-c; z)$$ also satisfies this HGDE, and if $c\neq 1$ then this function and $\F(a,b;c; \lambda)$ are a basis of solutions in any simply connected subset of the region of convergence $|\lambda|<1$. These functions may be analytically continued to a basis of solutions in a neighborhood of any point in $\CP\setminus\{0,1,\infty\}$.


\subsection{Schwarz triangle mappings.}\label{SS:triangle} Let $\bH\subset \CP\setminus\{0,1,\infty\}$ denote the upper half-plane, and let $f$ and $g$ be two linearly independent solutions to \ref{E:DE} defined on $\bH$. The Schwarz triangle map is defined as $$D:\bH \to \CP, \quad D(\lambda)=\frac{f(\lambda)}{g(\lambda)}.$$ Define
\begin{equation}\label{E:angles}
\kappa=|1-c|, \quad \mu=|c-a-b|, \quad \nu=|a-b|.
\end{equation}

\begin{thm}\label{T:triangle}
Assume  $0\leq \kappa, \mu, \nu<1$. The Schwarz triangle mapping $D$ maps $\bH$ diffeomorphically to a triangle in $\CP$ bounded by segments of circles or lines, with angles $\pi\kappa, \pi\mu$ and $\pi\nu$ at $D(0), D(1)$ and $D(\infty)$ respectively. 
\end{thm}

For the convenience of the reader, we now sketch the proof. Readers willing to accept Theorem \ref{T:triangle} as a black box are advised to proceed directly to the next section. 

The Wronskian is $W(f,g)=fg'-f'g$, and by Abel's identity it can be computed to be, up to a nonzero constant multiple, $\lambda^{-c}(\lambda-1)^{c-a-b-1}$. Since $D'(\lambda)=W(f,g)/g^2$, we see that the derivative of the Schwarz map is nonzero at all $\lambda$ where $g(\lambda)\neq 0$. At a point $\lambda_0\neq 0,1$ with $g(\lambda_0)=0$, we instead consider $1/D(\lambda)$; $f$ and $g$ may not simultaneously vanish away from the singular points $0,1, \infty$ (if they did, the Wronkskian would vanish there, and hence everywhere). We conclude that the map $D$ is a local diffeomorphism.

On any of the real arcs $(-\infty, 0), (0,1), (1, \infty)$, the HGDE \ref{E:DE} has a basis $f_0,g_0$ of real valued solutions, for example, derived from the above hypergeometric functions if $c\neq 1$. Being a real valued function, $D_0=\frac{f_0}{g_0}$ maps the chosen arc to a real arc (connected interval of $\bR$). Since the HGDE has a two dimensional space of solutions, there is a matrix $A\in GL_2(\bC)$ so that 
$$\left( \begin{array}{c} f \\ g\end{array}\right)= A \left( \begin{array}{c} f_0 \\ g_0\end{array}\right).$$
On the chosen arc we conclude that $D(\lambda)=A(D_0(\lambda))$, where here $A$ acts as a M\"{o}bius transformation. Hence $D$ maps the real arcs $(-\infty,0)$, $(0,1),$ and $(1,\infty)$ to segments of circles or straight lines.

If $c\neq 1$, set $$f_0=\lambda^{1-c}\F(a+1-c, b+1-c; 2-c; z)\quad\text{and}\quad g_0= \F(a,b;c; \lambda).$$ As before there is a $A\in GL_2(\bC)$ such that $D(\lambda)=A(D_0(\lambda))$. The map $D_0$ is, for $\lambda$ near $0$, approximately $\lambda^{1-c}$. Hence $D$ sends a neighborhood of $0$ in $\bH$ to a corner with angle $\kappa \pi$. If $c=1$, then $f_0$ must be replaced with a logarithmic solution of the form $\log(z)g_0$ plus a holomorphic function. In this case, $D_0$ is, for $\lambda$ near $0$, approximately $\log(z)$ plus a constant. Hence $D$ sends a neighborhood of $0$ in $\bH$ to a corner with angle $0=\kappa \pi$.

Similarly it can be shown that $D$ maps neighborhoods of $1$ and $\infty$ in $\bH$ to corners with angles $\mu\pi$ and $\nu\pi$. In summary, $D:\bH\to\CP$ is a local diffeomorphism which can be extended to the closure $\overline{\bH}$ of $\bH$ in $\CP$; this extension maps the boundary of $\bH$ bijectively to the boundary of the desired triangle, and maps a neighborhood in $\bH$ of the boundary (not including the boundary to itself) to such an interior neighborhood of the triangle. Any such map $D$ can be shown to be a diffeomorphism onto the triangle.


\section{Period mapping of $H^1(\br)$}\label{S:periodmap}

In this section we describe the period mapping for each subbundle $H^1(\br)$ of $H^1$.  The Gauss-Manin connection can be used to canonically identify nearby fibers of $H^1(\br)$. Here we will work over $\bH\subset \CP\setminus\{0,1,\infty\}$. Because this set is simply connected, we may assume that all fibers of $H^1(\br)$ have been thus identified. As $\lambda\in \CP\setminus\{0,1,\infty\}$ varies, the position of $H^{1,0}(\br)$ in $H^1(\br)$ varies. 

The period mapping is a holomorphic function to the Grassmanian of $\dim H^{1,0}(\br)$ planes in a $\dim H^1(\br)$ dimensional complex vector space which records this \emph{variation of Hodge structure}. If $H^1(\br)$ is equal to $H^{1,0}(\br)$ or $H^{0,1}(-\br)$ this map is trivial. If $\dim H^{1,0}=1=\dim H^{0,1}$, then the position of the one (complex) dimensional space $H^{1,0}$ inside the two dimensional space $H^1(\br)$ is recorded by it's \emph{slope}
$$\tau(\lambda)=\frac{ \int_{\alpha} \omega_\lambda}{ \int_{\beta}\omega_\lambda},$$
where $\omega_\lambda \in H^{1,0}$ and $\alpha$ and $\beta$ are two fixed homology classes. The map $\tau$ is also called the period mapping. Very often, it's image can be assumed to lie in the upper half plane $\bH$. 

We suggest \cite[Section 1.1]{C} as an accessible introduction to the relations between period mappings, monodromy, and differential equations.


\subsection{Sections satisfy a differential equation.} Fix the family $\M_N(A)$, and fix $\br$ in the row span of $A$. Assume $t(-\br)=2=t(\br)$, so $\dim H^{1,0}(\br)=1=\dim H^{0,1}(\br)$ (Proposition \ref{P:dim}). Set $t_i=t_i(-\br)$. The bundle $H^{1,0}(\br)$ has a global section given by Proposition \ref{P:dim} 
\begin{equation}\label{E:omega}
\omega=z^{-t_1}(z-1)^{-t_2}(z-\lambda)^{-t_3}dz.
\end{equation} We compute the differential equation satisfied by the section $\omega$. The first and second covariant derivatives of $\omega$ in the direction $\partial_\lambda$, with respect to the Gauss-Manin connection, are 
\begin{eqnarray*}
\omega' 
&=&
t_3 z^{-t_1}(z-1)^{-t_2}(z-\lambda)^{-t_3-1} dz 
\\&=&
 \frac{t_3}{z-\lambda}\omega, 
\\
\omega'' 
&=&
t_3 (t_3+1) z^{-t_1}(z-1)^{-t_2}(z-\lambda)^{-t_3-2} dz
\\&=&
 \frac{t_3+1}{z-\lambda}\omega' = \frac{t_3(t_3+1)}{(z-\lambda)^2}\omega.
\end{eqnarray*}

We wish to find a relation between $\omega, \omega'$  and $\omega''$ in first cohomology. These three cohomology classes all lie in the two dimensional space $H^1(\br)$, so such a relation necessarily exists. We compute: 
\begin{eqnarray*}
0 
&=&
d\left[ -t_3 z^{-t_1+1}(z-1)^{-t_2+1}(z-\lambda)^{-t_3-1} \right]
 \\&=&
 \lambda(\lambda-1)\omega'' + \left[ (t_1+t_2+2t_3)\lambda-(t_1+t_3)\right]\omega'
 \\&&
 +t_3(t_1+t_2+t_3-1)\omega.
\end{eqnarray*}

Hence, the differential equation satisfied by the global section $\omega$ is 
\begin{equation}\label{E:GM}
\lambda(\lambda-1)\nabla^2+\left[(a+b+1)\lambda-c\right]\nabla + ab = 0,
\end{equation}
where $\nabla$ is the Gauss-Manin connection, and 
\begin{equation}\label{E:abc}
a=t_1+t_2+t_3-1, \quad b=t_3, \quad c=t_1+t_3.
\end{equation}

Note that our assumption $t(\br)+t(-\br)=4$ implies that all $t_j\neq 0$. (In general, $t(\br)+t(-\br)$ is equal to the number of nonzero $t_j$.) It follows that $c>0$, and $c$ is an integer if and only if $c=1$.  

Compute 
\begin{eqnarray} \label{E:angles2}
\kappa&=&|1-c|=|1-t_1-t_3|,\\
\mu&=&|c-a-b|=|1-t_2-t_3|,\\
\nu&=&|a-b|=|1-t_1-t_2|.
\end{eqnarray}
Since the $t_j$ are nonzero we have $0\leq \kappa, \mu, \nu <1$.


\subsection{Period Mappings.} Real homology classes $\alpha$ and $\beta$ can be chosen so that the functions $g_\alpha(\lambda)= \int_{\alpha} \omega$ and $g_\beta(\omega)= \int_{\beta}\omega$ are nonconstant and noncollinear. Consider the period map 
$$\tau(\lambda)=\frac{ \int_{\alpha} \omega}{ \int_{\beta}\omega}.$$ 
 All of these functions are holomorphic. Since the section $\omega$ satisfies \ref{E:GM}, the functions $g_\alpha$ and $g_\beta$ must satisfy the corresponding HGDE \ref{E:DE}. (The only difference between \ref{E:DE} and \ref{E:GM} is the substitution of the Gauss-Manin connection.) It follows that the period mapping $\tau$ is a Schwarz triangle map.
 
\begin{prop}\label{P:pmap}
Let $\br\neq0 $ be in the row span of $A$, and assume $t(-\br)=2=t(\br)$. Then the period map $\tau$ restricted to $\bH\subset \CP\setminus\{0,1,\infty\}$ may be assumed to map biholomorphically onto a hyperbolic triangle with angles $\pi\kappa, \pi\mu, \pi\nu$ respectively at $\tau(0), \tau(1), \tau(\infty)$. Here
$$\kappa=|1-t_1-t_3|,\quad
\mu=|1-t_2-t_3|,\quad
\nu=|1-t_1-t_2|.$$
\end{prop}

\begin{proof}
An elementary computation gives that $\kappa+\mu+\nu<1$ when $t(-\br)=2=t(\br)$ (see the proof of Theorem \ref{T:arealyaps}). Since $\tau$ is a Schwarz triangle mapping, it maps $\bH$ to triangle with angles $\pi\kappa, \pi\mu, \pi\nu$. By applying a M\"obius transformation, this triangle may be assumed to lie in $\bH$. Applying this M\"obius transformation corresponds to picking a different choice of $\alpha$ and $\beta$ in the definition of $\tau$ above. 
\end{proof}


\section{Lyapunov exponents}\label{S:lyaps} 

In this section we prove Theorems \ref{T:periodlyap} and \ref{T:arealyap}, and give an algorithm for computing the Lyapunov exponents of the bundle $H^1$ for abelian square-tiled surfaces.


\subsection{Background.} We recommend the introduction of \cite{EKZbig} for a survey on Lyapunov exponents of flat surfaces. Here we will restrict ourselves to the case where we have an equivariant (flat) rank 2 subbundle $W\subset H^1$ with $W^{1,0}=\overline{W^{0,1}}$ having rank 1. Furthermore, we will assume that the base $\M$ of this bundle is a finite cover of a Teichm\"uller curve. 

Note that, unlike in previous treatments, we use complex cohomology instead of real cohomology, and the subbundle $W$ is a complex subbundle of $H^1$ which need not come from a real subbundle of real cohomology. This difference is of no consequence in the results that we quote.

Kontsevich, Zorich, and Forni \cite{F} have shown the following formula for the nonnegative Lyapunov exponent $\ell$ of $W$.
\begin{equation}\label{E:lyap}
\ell = \frac{-1}{\operatorname{Area}(\M)} \int_{\M} \Delta_{hyp} \log \left| \langle \omega_\lambda, \omega_\lambda  \rangle \right| dg_{hyp}(\lambda)
\end{equation}
Since $\M$ is a finite cover of a Teichm\"uller curve, it is a hyperbolic orbifold, whence the reference to its area, the hyperbolic Laplacian $\Delta_{hyp}$, and the hyperbolic area form $dg_{hyp}$. In the above formula, $\omega_{\lambda}$ is a global section of $H^{1,0}$; if none exists, the contributions of local sections can be added up using a partition of unity. The form $\langle \cdot, \cdot \rangle$ is the usual hermitian intersection form.  

In fact, there is a formula for $\ell$ in terms of the degree of the Deligne extension of $W^{1,0}$.  Formula \ref{E:lyap} is an intermediate step in the proof of this formula. For a proof see the Teichm\"uller curve case of \cite{EKZbig}, or see \cite[Section 9]{BM}.

Note that $\left| \langle \omega_\lambda, \omega_\lambda \rangle\right|$ is the Hodge norm of $\omega_\lambda \in H^{1,0}$. 

Since the Kontsevich-Zorich cocycle preserves a symplectic form (the usual intersection form), its Lyapunov spectrum is symmetric under negation.


\subsection{Lyapunov exponents from period mappings.} Let $W$ be as above. The universal cover of the Teichm\"uller curve $\M$ is the upper half-plane $\bH$. Let $\omega_\lambda$ be a global section of the pull back of $W^{1,0}$ to $\bH$. The subscript indicates the point $\lambda\in \bH$. 

Let $\alpha^*, \beta^*$ be a pair of global flat sections for the pull back of $H^1$ to  $\bH$, with intersection pairings 
\begin{equation}\label{E:pairs}
\langle \alpha^*, \beta^*\rangle = -\langle \beta^*, \alpha^* \rangle=1 \quad \text{      and      }\quad \langle \alpha^*, \alpha^*\rangle = \langle \beta^*, \beta^*\rangle=0.
\end{equation}
Here we use the dual of the usual hermitian intersection form on $H_1$.

Picking $\alpha^*$ and $\beta^*$ appropriately, we may assume that the period map $\tau: \bH\to \CP$ $$f(\lambda) = \frac{\langle \omega_\lambda, \alpha^*\rangle}{\langle \omega_\lambda, \beta^* \rangle}$$ has image in $\bH$. Under this assumption, $\tau $ is well defined up to postcomposition by hyperbolic isometries. Hence the hyperbolic norm of its derivative,
$$\| \tau'(\lambda)\|_{hyp}= \frac{\Im(\lambda)}{\Im(\tau(\lambda))} | \tau'(\lambda)|,$$
is in fact a well defined function on $\M$.

\begin{proof}[\emph{\textbf{Proof of Theorem \ref{T:arealyap}}}]
In the definition of $\tau $, we may pick $\alpha^*, \beta^*$ to span $W$ (and so that \ref{E:pairs} holds). It then follows from the definition of $\tau $ that we may take
$$\omega_\lambda= \alpha^*-\tau(\lambda)\beta^*.$$ 
Let $F$ be a fundamental domain for $\M$ in $\bH$. By formula \ref{E:lyap}, the nonnegative Lyapunov exponent of $W$ is 
$$\ell=\frac{-1}{\operatorname{Area}(F)} \int_{F} \Delta_{hyp} \log \left| \langle \omega_\lambda, \omega_\lambda  \rangle \right| dg_{hyp}.$$

Write $\lambda=x+iy$, and $\tau =u+iv$, and compute that $$|\langle \omega_\lambda, \omega_\lambda  \rangle|=2v.$$

Recall $\Delta_{hyp}=y^2(\partial_x^2+\partial_y^2)$.  The Cauchy-Riemann equations give $$-\Delta_{hyp} \log |\langle \omega_\lambda, \omega_\lambda \rangle|= \frac{y^2}{v^2} |\tau'|^2,$$ which establishes the result. 
\end{proof}


\subsection{Lyapunov exponents of $\M_N(A)$.}\label{SS:lyaps} We return to the situation of an abelian square-tiled surface $M_N(A)$ and the corresponding finite cover $\M=\M_N(A)$ of a Teichm\"uller curve. The following is an expanded version of Theorem \ref{T:arealyap}. 

\begin{thm}\label{T:arealyaps}
Let $\br$ be in the row span of $A$. Unless $t(-\br)=t(\br)=2$, the Lyapunov exponents of $H^1(\br)$ are zero. So suppose $t(-\br)=2=t(
\br)$, and the period map restricted to $\bH\subset \CP\setminus\{0,1,\infty\}$ is a biholomorphism onto a triangle with hyperbolic area $A$. In this case the Lyapunov spectrum of $H^1(\br)$ is $\pm\frac{A}{\pi}$. Furthermore
$$\frac{A}{\pi}=2\min_{j=1,2,3,4} \{t_j(-\br), 1-t_j(-\br)\}.$$
\end{thm} 

\begin{proof}
The dimension of $H^1(\br)$ is $t(-\br)+t(\br)-2$, and the dimension of $H^{1,0}(\br)$ is $t(-\br)-1$ (Proposition \ref{P:dim}). Unless $t(-\br)=t(\br)=2$, the bundle $H^1(\br)$ is equal to its $(1,0)$ or $(0,1)$ part and the restriction of the Hodge bilinear form \ref{E:hodge} to $H^1(\br)$ is definite. It follows that the monodromy of the Gauss-Manin connection lies in a compact group and the Lyapunov exponents are zero. 

In the remaining case case, $t(\br)=t(-\br)=2$, so $\dim H^{1,0}(\br)=\dim H^{0,1}(\br)=1$. 

The (finite cover of a) Teichm\"uller curve $\M$ is isometric to the hyperbolic three times punctured sphere, which we write as $\CP\setminus\{0,1,\infty\}$ (remark \ref{R:3ps}). Moreover, the period map of $H^1(\br)$ is a Schwarz triangle map $f$ which maps $\bH\subset \CP\setminus\{0,1,\infty\}$ to a triangle of area $A=1-\kappa-\mu-\nu$, where the angles $\kappa, \mu, \nu$ are given in Proposition \ref{P:pmap}. 

The hyperbolic three times punctured sphere is composed to two ideal triangles glued along edges. The two ideal triangles can be taken to be $\bH, \overline{\bH}\subset \CP\setminus\{0,1,\infty\}$.

Theorem \ref{T:periodlyap} gives $$\ell=\frac1{Area(\M)} \int_\M \|f'(\lambda)\|^2_{hyp} dg_{hyp}.$$
Instead of integrating over $\M$, we may integrate over a fundamental domain $F\subset \bH$ consisting of two ideal triangles $T_1$ and $T_2$. The image of each of these triangles (each of area $\pi$) is a triangle of area $A$ (figure \ref{F:triangleslyap}). Furthermore, the restriction of $f$ to each $T_i$ is a diffeomorphism onto its image.  The change of variables formula gives that $\int_{T_i} \|f'(\lambda)\|^2_{hyp} dg_{hyp}$ is the hyperbolic area $A$ of $f(T_j)$, so $\ell= \frac{A}{\pi}$. 

\begin{figure}[h]
\includegraphics[scale=0.5]{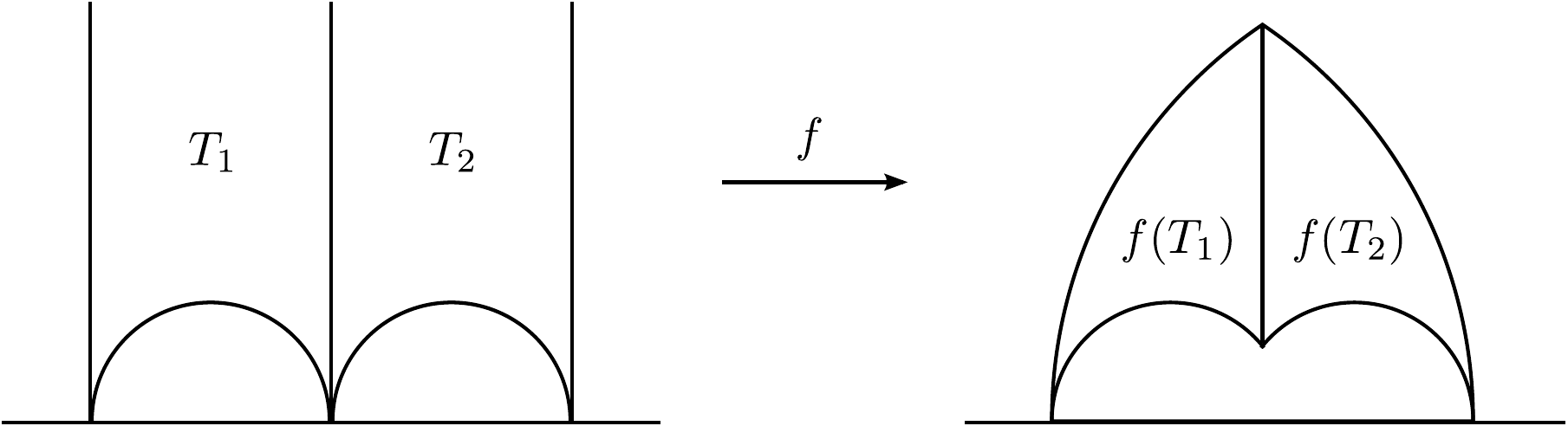}
\caption{The period map $f$ maps a fundamental domain $F=T_1\cup T_2$ biholomorphically to two triangles, each of area $A$.}
\label{F:triangleslyap}
\end{figure}

The restriction of the symplectic form coming from the intersection form to $H^1(\br)$ is symplectic, so the Lyapunov spectrum of $H^1(\br)$ is symmetric under negation.

Set $t_j=t_j(-\br)$.
$$\frac{A}{\pi}=1-(\kappa+\mu+\nu)=1- |1-t_1-t_3|-|1-t_2-t_3|-|1-t_1-t_2|.$$
Under the assumption that $t(\bk)=t(-\bk)=2$, this is symmetric in the $t_j$, so we may assume without loss of generality that $t_1\leq t_2\leq t_3\leq t_4$. There are now two cases: the first is $t_1\leq 1-t_4$. In this case, the expression for $\lambda=\frac{A}{\pi}$ is $$1-(1-t_1-t_3)+(1-t_2-t_3)-(1-t_1-t_2)= 2t_1.$$
In the second case, $t_1> 1-t_4$, and the expression for $\frac{A}{\pi}$ simplifies to $2(1-t_4)$. 
\end{proof}

Here is an algorithmic restatement of this theorem, which has been included for comparison with cyclic case in \cite{EKZsmall}. 

\begin{thm}\label{T:alglyap}
Let $M_N(A)$ be an abelian square-tiled surface. Start with $\Lambda=\emptyset$, and for every $\br$ in the row span of $A$, 
\begin{itemize}
\item if $t(-\br)=2=t(\br)$, add $2\min_{j=1,2,3,4} \{t_j(-\br), 1-t_j(-\br)\} $ to $\Lambda$,
\item if $(-\br)=3$, add $0$ to $\Lambda$. 
\end{itemize} 
The resulting set $\Lambda$ is the nonnegative part of Lyapunov spectrum of $H^1$.
\end{thm}



\bibliography{mybib}{}
\bibliographystyle{amsalpha}
\end{document}